\documentclass[final,3p,sort&compress]{elsarticle}

\usepackage{epsfig}
\usepackage{amssymb}
\usepackage{amsmath}
\usepackage{amsthm}
\usepackage{amsbsy}
\usepackage{multirow}
\usepackage{graphicx}
\usepackage[pdftex,bookmarksnumbered]{hyperref}
\usepackage{multicol}

\usepackage{setspace}
\usepackage{color}

\usepackage{lineno}

\begin{document}

\begin{frontmatter}



\title{Strict upper and lower bounds for quantities of interest in static response sensitivity analysis}


\author{Mengwu Guo\corref{cor1}}
\ead{gmw13@mails.tsinghua.edu.cn}
\cortext[cor1]{Corresponding author.}

\author{Hongzhi Zhong\corref{}}
\ead{hzz@tsinghua.edu.cn}

\address{Department of Civil Engineering, Tsinghua University, Beijing 100084, China}

\begin{abstract}
   In this paper, a goal-oriented error estimation technique for static response sensitivity analysis is proposed based on the constitutive relation error (CRE) estimation for finite element analysis (FEA). Strict upper and lower bounds of various quantities of interest (QoI) that are associated with the response sensitivity derivative fields are acquired. Numerical results are presented to assess the strict bounding properties of the proposed technique.

\end{abstract}




\end{frontmatter}


\section{Introduction}

In the design of engineering structures, the finite element method (FEM) has been widely used to make critical decisions. In order to control the quality of numerical simulations and develop confidence in decisions, a research topic, referred to as
\emph{model verification}, has been intensively studied for more than four decades. Among different error sources of numerical simulations for a chosen model, the discretization error is predominant and controllable. For the purpose of evaluating discretization error in finite element analysis (FEA), several families of \emph{a posteriori} error estimators \cite{error1,error2,error3,error4} have been presented for the estimation of errors measured in global norms, such as explicit error estimators \cite{posterr0.5}, implicit error estimators \cite{posterr1,posterr2}, recovery-based error estimators \cite{posterr3}, hierarchical estimators \cite{posterr4}, \emph{constitutive relation error} (CRE) estimators \cite{posterr5}, etc.

The goal of many finite element computations in structural analysis is the determination of some specific quantities of interest, such as local stress values, displacements etc., which is necessary for a particular design decision. Thus, it is frequently the case that \emph{a posteriori} finite element error analysis is focused on goal-oriented error estimation. Towards this end, adjoint/dual-based techniques are used to estimate the errors in solution outputs, which have been systematically reviewed in \cite{adj1,adj2,adj3,adj4}. Research on goal-oriented error estimation was initiated in the 1990's \cite{goalerr0,goalerr1,goalerr2,goalerr3,goalerr4,goalerr5,goalerr5.1,goalerr5.2}. Since then, several methods have been developed and applied to solutions of various problems, such as Poisson's equation, linear and non-linear static problems in solid mechanics, eigenvalue problems, time-dependent problems, non-trivial problems of CFD, etc (see \cite{goalerr9,goalerr10} for example). A variety of specific error estimation techniques have been proposed to evaluate the discretization error in quantities of interest, for instance, the adjoint-weighted residual method \cite{adj1,adj4,goalerr9}, the energy norm based estimates \cite{goalerr6}, the Green's function decomposition method \cite{goalerr7}, the strict-bounding approach based on Lagrangian formulation \cite{goalerr8}, the CRE-based error estimation \cite{goalerr5}. Among the available techniques, the CRE-based error estimation provides guaranteed strict bounds of quantities. The strict bounding property, together with its advantage of wide applicability \cite{CRE1,CRE2,CRE3,CRE4,CRE5,CRE6,CRE7,CRE8,CRE9,Guo1,Guo2}, makes the CRE stand out for goal-oriented error estimation.

Sensitivity analysis plays an important role in uncertainty analysis, structural optimization and many other areas of structural analysis. When using some numerical methods in the first-order perturbed formulation to compute the static response sensitivity of a structural system, discretization error exists in the analysis. For instance, the stochastic perturbation method \cite{pert} is usually chosen to obtain statistically characteristic values of some structural outputs, in which the response sensitivity derivatives with respect to input parameters appear in the expressions of coefficients. Hence controlling the discretization error in response sensitivity helps enhance the accuracy of evaluating the statistically characteristic outputs. In the context of structural optimization or other parameterized problems that require repeatedly solving the structural responses under different inputs, gradient-based algorithms desire the response sensitivity derivatives at each iteration step in the parameter space. If some reduced order methods, such as the reduced basis method \cite{RB1} and the proper generalized decomposition \cite{PGD1}, are used to solve the structural responses and response derivatives at a number of sampling points with a decreased computational cost, the verification of numerical simulations will also play a crucial part throughout the procedure, see \cite{RB1,RB2,PGD2} for examples.
Therefore, \emph{a posteriori} estimators are required to estimate the discretization error in the solution for sensitivity derivatives of the structural response, especially in some specific quantities about the response sensitivity. As far as the authors know, the relevant error estimation techniques have not been adequately studied, and only limited information has been available. For example, an explicit (residual-based) error estimator has been used in \emph{a posteriori} error estimation in sensitivity analysis \cite{sen1}, and a Neumann-subproblem \emph{a posteriori} finite element procedure has been proposed to provide upper and lower bounds for functionals of the response sensitivity derivative fields \cite{sen2}.

On the basis of the principle of minimum complementary energy, the CRE-based goal-oriented error estimation will be extended to the cases of non-symmetric bilinear forms, especially to the static response sensitivity analysis of linear structural systems by the FEM in this paper. Consequently, strict upper and lower bounds can be obtained for quantities of interest, which are linear functionals associated with the sensitivity derivative fields of displacements, including the sensitivity derivatives of some scalar-valued static response quantities.

Following the introduction, the basics of the CRE estimation and the CRE-based goal-oriented error estimation are reviewed in Section 2. In Section 3, the CRE-based error estimator is extended to the cases with non-symmetric bilinear forms, and in Section 4, the estimator is used for goal-oriented error estimation of static response sensitivity. Numerical results for some model problems are presented to assess bounding property of the proposed estimation technique in Section 5. In Section 6, conclusions are drawn.

\section{Basics of the constitutive relation error estimation}

%
%

\subsection{An abstract primal problem}

To start with, a typical problem in structural analysis is introduced \cite{abstract}. A Banach space $\mathcal{V}$, referred to as the 'space of kinematically admissible solutions', consists of all the possible displacements that satisfy the Dirichlet boundary conditions \footnote{In this paper, only the problems with homogeneous Dirichlet boundary conditions are discussed, since those with nonhomogeneous Dirichlet boundary conditions can be equivalently transformed to homogeneous cases.}. As its dual space, the 'loading space' $\mathcal{V}^*$ is given with the duality pair $_{\mathcal{V}^*}\langle\cdot,\cdot\rangle_{\mathcal{V}}$. Usually, a load $f\in \mathcal{V}^*$ includes a body force in the domain that the structure occupies and a traction on its Neumann boundary. A Banach space of strains, $\mathcal{E}$, and its dual space~--~the space of stresses, $\mathcal{E}^*$, are introduced, and their duality pair is written as $_{\mathcal{E}^*}\langle\cdot,\cdot\rangle_{\mathcal{E}}$.

The relation between a displacement element $v\in \mathcal{V}$ and its corresponding strain $\varepsilon \in \mathcal{E}$ is represented by a linear differential operator $A:\mathcal{V}\to\mathcal{E}, v\mapsto \varepsilon$, i.e. $\varepsilon=Av$. The adjoint operator of  $A$,  denoted by $A^*:\mathcal{E}^*\to \mathcal{V}^*$, is then defined as
\begin{equation}\label{adjoint}
_{\mathcal{E}^*}\langle \tau,Av\rangle_{\mathcal{E}}=~_{\mathcal{V}^*}\langle A^*\tau,v\rangle_{\mathcal{V}}\,\quad \forall(\tau,v)\in \mathcal{E}^*\times \mathcal{V}\,.
\end{equation}
In structural analysis, the operator $A$ is usually gradient-like and $A^*$ is divergence-like, which is a natural derivation from Green's formula. Besides, the relation between stresses and strains, or termed the constitutive relation, is represented by a material operator $K:\mathcal{E}\to \mathcal{E}^*$.

Then the governing equations for the primal structural problem are given as follows:
\begin{equation}
\begin{split}
& u\in \mathcal{V}\,,\;\varepsilon\in\mathcal{E}\,,\;\sigma\in\mathcal{E}^*\,,\\
& \varepsilon=Au\,,\;\sigma=K\varepsilon\,,\;A^*\sigma=f\,,
\end{split}
\end{equation}
or written with a single unknown $u$ as
\begin{equation}\label{gov}
u\in \mathcal{V}\,,\;A^*(K(Au))=f\,.
\end{equation}
With the aid of Eq.~\eqref{adjoint}, the weak form of Eq.~\eqref{gov} is stated as: find $u\in \mathcal{V}$ such that
\begin{equation}\label{weak}
_{\mathcal{E}^*}\langle K(Au),Av \rangle_{\mathcal{E}}=~_{\mathcal{V}^*}\langle f,v \rangle_{\mathcal{V}}\,\quad \forall v\in \mathcal{V}\,,
\end{equation}
which is also referred to as the virtual work principle. 


In this paper, linear elastic problems are taken into consideration. Then $\mathcal{V}$, $\mathcal{V}^*$ and $\mathcal{E}=\mathcal{E}^*$ are ascribed to Hilbert spaces, $_{\mathcal{E}^*}\langle \cdot,\cdot \rangle_{\mathcal{E}}$ is symmetric and positive definite, and the operator $K$ is linear, reversible, symmetric and positive definite. In this case, $A^*KA$, the differential operator in Eq.~\eqref{gov} is of elliptic type. For example, $A=A^*=\partial_{xx}$ and $K=EI(x)$ for a beam problem, where $EI$ is the flexible stiffness; $A=(\nabla+\nabla^T)/2$, $A^*=-\mathrm{div}$ and $K$ is the Hooke's stiffness tensor for a 2D or 3D problem in linear elasticity.

For notation, a symmetric semi-positive definite bilinear form $a_{{u}}(\cdot,\cdot): \mathcal{V} \times \mathcal{V} \to \mathbb{R}$ and the corresponding semi-norm, a symmetric positive definite bilinear form $a_{{\sigma}}(\cdot,\cdot): \mathcal{E}^* \times \mathcal{E}^* \to \mathbb{R}$ and the corresponding norm are introduced, respectively, as
\begin{equation}
\begin{split}
&  a_{{u}}(u,v)=_{\mathcal{E}^*}\langle K(Au),Av \rangle_{\mathcal{E}}\,,\;\|{u}\|_{{u}}= \sqrt{a_{{u}}(u,u)}\,,\;(u,v)\in \mathcal{V}\times\mathcal{V};\\
&  a_{{\sigma}}({\sigma},{\tau})=_{\mathcal{E}^*}\langle \sigma,K^{-1}\tau \rangle_{\mathcal{E}}
\,,\;\|{\sigma}\|_{{\sigma}}= \sqrt{a_{{\sigma}}({\sigma},{\sigma})}\,,\; (\sigma,\tau)\in \mathcal{E}^*\times \mathcal{E}^*\,.
\end{split}
\end{equation}
The duality pair $_{\mathcal{V}^*}\langle\cdot,\cdot\rangle_{\mathcal{V}}$ is then written as $\langle\cdot,\cdot\rangle$ for simplification. It can be proven that $a_u$ is a continuous and coercive bilinear form for linear elastic problems. This ensures the existence and uniqueness of the solution to Eq.~\eqref{weak}, which is restated as: find $u\in\mathcal{V}$ such that
\begin{equation}
a_u(u,v)=\langle f,v\rangle\quad \forall v\in\mathcal{V}\,.
\end{equation}

The primal problem \eqref{weak} can be formulated as follows: find a displacement field $u$ and a stress field ${\sigma}$ satisfying
\begin{itemize}
  \item The compatibility condition:
   \begin{equation}\label{k}
   u \in \mathcal{V}\,;
   \end{equation}
  \item The equilibrium condition:
   \begin{equation}\label{e}
    {\sigma} \in \mathcal{E}^*,\;  a_{{\sigma}}({\sigma},K(A(v)))=\langle f,v\rangle\quad\forall {v} \in \mathcal{V}\,;
   \end{equation}
  \item The constitutive relation: (Hooke's law)
   \begin{equation}\label{c}
   {\sigma}=K(A(u))\,.
   \end{equation}
\end{itemize}

When Eqs. \eqref{k}, \eqref{e} and \eqref{c} hold true, the pair $({u},{\sigma})$ is the exact solution to the primal problem. To seek numerical solutions, the problem can be discretized using the displacement-based Galerkin finite element method, i.e. find ${u}_{h} \in \mathcal{V}_{h}$ such that
\begin{equation}\label{fe}
  a_{{u}}({u}_{h},{v})= \langle f,v\rangle\quad\forall {v} \in \mathcal{V}_{h}\,,
\end{equation}
where $ \mathcal{V}_{h}=\mathcal{V}\cap \mathcal{P}_{h}$, and $\mathcal{P}_{h}$ denotes the finite element space under the mesh characterized by size $h$. Together with the solution of stress field
\begin{equation}
{\sigma}_{h}=K(A(u_h))\,,
\end{equation}
in the sense of distribution, the pair $({u}_{h},{\sigma}_{h})$ forms the finite element approximations of the primal problem, resulting in a discretization error.

\subsection{Concept of constitutive relation error}

The concept of constitutive relation error (CRE) relies on the concept of
admissible solution pair $(\hat{{u}},\hat{{\sigma}})$, i.e. the combination of the kinematically admissible field $\hat{{u}}$ verifying \eqref{k} and
the statically admissible field $\hat{{\sigma}}$ verifying \eqref{e}.
The solution quality is quantified by the error of constitutive relations.

Hence an error measured in terms of the constitutive relation is defined as
\begin{equation}
   e_{CRE}(\hat{{u}},\hat{{\sigma}}):=\| \hat{{\sigma}}- K(A(\hat{{u}}))\|_{{\sigma}},
\end{equation}
which is the \textit{constitutive relation error} (CRE).

An important property of the constitutive relation error is the Prager-Synge theorem \cite{PS}:
\begin{equation}\label{bdg}
e_{CRE}^{2}(\hat{{u}},\hat{{\sigma}})=\|\hat{{\sigma}}- {\sigma}   \|_{{\sigma}}^{2}+\| \hat{{u}}-{u}\|_{{u}}^{2}\,.
\end{equation}
Then a corollary that $e_{CRE}(\hat{{u}},\hat{{\sigma}})=0 \Leftrightarrow (\hat{{u}},\hat{{\sigma}})=({u},{\sigma})$ a.e. follows immediately.

\subsection{Global discretization error estimation based on the CRE}

The finite element solution for displacements satisfies ${u}_{h} \in \mathcal{V}_{h} \subset \mathcal{V}$, meaning that ${u}_{h}$ can be taken as the kinematically admissible field, i.e. $\hat{{u}}={u}_{h}$. However, the finite element solution for stresses ${\sigma}_{h}$ does not satisfy the equilibrium equations, i.e. $\hat{{\sigma}}\neq{\sigma}_{h}$. There already exist plenty of techniques proposed to recover the equilibrated stress field  $\hat{{\sigma}}=\hat{{\sigma}}_{h}$ from the finite element stress solution $\sigma_h$ via an energy relation called \textit{prolongation condition}, see \cite{error1,CRE10} for reviews.

According to Eq. \eqref{bdg}, the constitutive relation error $e_{CRE}({u}_{h},\hat{{\sigma}}_{h})$, which can be considered as a global discretization error estimator, provides an upper bound of the global energy norm error of the finite element solution, i.e.
\begin{equation}\label{ubp}
\|{u}-{u}_{h}\|_{{u}}=\|{\sigma}-{\sigma}_{h}\|_{{\sigma}}\leq e_{CRE}({u}_{h},\hat{{\sigma}}_{h}).
\end{equation}

As a matter of fact, this bounding property \eqref{ubp} is guaranteed by the well-known principle of minimum complementary energy. Notice that ${e}:={u}-{u}_h$ is the solution of such a 'residual' problem: find ${e}\in \mathcal{V}$ such that
\begin{equation}\label{errweak}
a_u(e,v)=\langle R,v\rangle\quad \forall {v}\in \mathcal{V}\,,
\end{equation}
where $R\in \mathcal{V}^*$ is defined as $\langle R,v\rangle=\langle f,v\rangle-a_{u}(u_h, v)$, ${v} \in \mathcal{V}$.
Then minimizing the complementary energy of this problem, also referred to as the dual variational formulation, immediately gives
the inequality \eqref{ubp}.

%
%

\subsection{Goal-oriented error estimation based on the CRE}

Assume that the quantity of interest is a linear bounded functional with respect to the displacement field ${u}$ defined in the global form $Q(u)=\langle Q,u\rangle$, where $Q \in \mathcal{V}^*$. Then, an adjoint problem associated with the output $Q(u)$ can be defined as: find $\tilde{{u}}\in \mathcal{V}$ such that
\begin{equation}
a_u(v,\tilde{u})=Q(v) \quad \forall v\in\mathcal{V}\,,
\end{equation}
or formulated as: find a displacement field $\tilde{u}$ and a stress field $\tilde{{\sigma}}$ that satisfy
\begin{itemize}
  \item The compatibility condition:
   \begin{equation}
   \tilde{{u}} \in \mathcal{V}\,;
   \end{equation}
  \item The equilibrium condition:
   \begin{equation}
   \tilde{{\sigma}} \in \mathcal{E}^*,\;  a_{{\sigma}}(K(A(v)),\tilde{\sigma})=Q(v)\quad \forall {v} \in \mathcal{V}\,;
   \end{equation}
  \item The constitutive relation:
   \begin{equation}
   \tilde{{\sigma}}=K^T (A(\tilde{u}))=K(A(\tilde{u}))\,.
   \end{equation}
\end{itemize}

Similar to the primal problem, the displacement field for the adjoint problem can be approximately obtained using the Galerkin finite element method: find $\tilde{{u}}_{h} \in \mathcal{V}_h$ such that
\begin{equation}
  a_u(v,\tilde{u}_h)=Q(v)\quad  \forall v \in \mathcal{V}\,,
\end{equation}
and the stress field solution is accordingly given as $\tilde{\sigma}_h=K(A(\tilde{u}_h))$ in the sense of distribution. Furthermore, an admissible pair $(\tilde{{u}}_{h},\hat{\tilde{{\sigma}}}_{h})$ for the adjoint problem can be derived using the same technique as that for the primal problem.

The approximation of quantity $Q(u)$ is usually computed as $Q(u_h)$.
With the fact that $u-u_h\in\mathcal{V}$ and $a_u(u-u_h,\tilde{u}_h)=0$, the error in quantity $Q(u)$ is given as
\begin{equation}
\begin{split}
Q(u)-Q(u_h)  = & ~Q(u-u_h)= a_{{u}}(u-u_h,\tilde{{u}}-\tilde{{u}}_h)\\
 =&  ~ \frac{1}{4}\left\|{\kappa}({u}-{u}_h)+\frac{1}{{\kappa}}(\tilde{{u}}-\tilde{{u}}_h)\right\|_{{u}}^{2}
 -\frac{1}{4}\left\|{\kappa}({u}-{u}_h)-\frac{1}{{\kappa}}(\tilde{{u}}-\tilde{{u}}_h)\right\|_{{u}}^{2},
\end{split}
\end{equation}
where $\kappa \in \mathbb{R}^{+}$ is an arbitrary parameter, and the parallelogram identity \cite{goalerr5.1,goalerr5.2} is used. Thus strict upper and lower bounds of $Q(u)-Q(u_h)$ can be represented in a computable form by introducing the admissible fields:
\begin{equation}\label{ded1}
\pm(Q(u)-Q(u_h))  \leq \frac{1}{4}\left\|{\kappa}({u}-{u}_h)\pm \frac{1}{{\kappa}}(\tilde{{u}}-\tilde{{u}}_h)\right\|_{{u}}^{2}
 \leq \frac{1}{4}\left\|{\kappa}(\hat{{\sigma}}_{h}-K(A(u_h)))\pm \frac{1}{{\kappa}}(\hat{\tilde{{\sigma}}}_{h}-K(A(\tilde{u}_h)))\right\|_{{\sigma}}^{2}\,.
\end{equation}
Taking $\kappa=\sqrt{\frac{e_{CRE}(\tilde{{u}}_{h},\hat{\tilde{{\sigma}}}_{h})}{e_{CRE}({u}_{h},\hat{{\sigma}}_{h})}}$, a pair of computable strict error bounds with the sharpest gap is given as follows:
\begin{equation}\label{Ibounds}
 \left|Q(u)-Q(u_h)-\frac{1}{2}a_{{\sigma}}(\hat{{\sigma}}_{h}-K(A(u_h)),\hat{\tilde{{\sigma}}}_{h}-K(A(\tilde{u}_h)))  \right|
\leq ~ \frac{1}{2}e_{CRE}({u}_{h},\hat{{\sigma}}_{h})\cdot e_{CRE}(\tilde{{u}}_{h},\hat{\tilde{{\sigma}}}_{h})\,.
\end{equation}

\section{Extension to cases of non-symmetric bilinear forms }\label{nonsym}

On the basis of the idea of splitting the operator into symmetric and antisymmetric parts, some output-based \emph{a posteriori} error bounds were proposed to deal with the problems with non-symmetric bilinear forms, such as the advection-diffusion-reaction problem \cite{nonsym1,nonsym2}. In this section, the symmetric part of a bilinear form is used to define \emph{the extended CRE-based goal-oriented error estimator of the problems with non-symmetric bilinear forms, which makes it possible to estimate the errors in quantities in static response sensitivity analysis.
}

The variational problem is usually stated as: find $u \in X$ such that
\begin{equation}
b_u(u,v)=l(v)\quad \forall v \in X,
\end{equation}
where $X$ is a Hilbert space, $b_u$ is continuous and coercive (not necessarily symmetric) bilinear form defined on $X \times X$, and $l$ a linear bounded functional on $X$ i.e. $l \in X^*$. In a finite element space $X_h \subset X$, an approximate solution $u_h$ can be found as
\begin{equation}
b_u(u_h,v)=l(v),\quad \forall v \in X_h.
\end{equation}

For a quantity of interest $Q(u)$ with $Q\in X^*$, the corresponding adjoint problem is defined as: find $\tilde{u}\in X$ such that
\begin{equation}
b_u(v,\tilde{u})=Q(v),\quad v\in X,
\end{equation}
and the corresponding finite element solution is denoted as $\tilde{u}_h \in X_h$.

Let us denote the symmetric part of $b_u$ as $b_u^{S}$, i.e. $b_u^S(u,v)=\frac{1}{2}(b_u(u,v)+b_u(v,u))$, $(u,v)\in X\times X$, and define $R_p \in X^*$ as $R_p(v)=l(v)-b_u(u_h,v)$ and $R_d \in X^*$ as $R_d(v)=Q(v)-b_u(v,\tilde{u}_h)$, $v \in X$. Since $b_u(u-u_h,u-u_h)=R_p(u-u_h)$ and $b_u(u-u_h,\tilde{u}_h)=0$, the error in the quantity $Q(u)$, i.e. $Q(u)-Q(u_h)=Q(u-u_h)$ can be represented as
\begin{equation}
\pm (Q(u)-Q(u_h))=\pm R_d(u-u_h)=\Pi_\kappa^\pm(\kappa (u-u_h)),
\end{equation}
where $\kappa \in \mathbb{R}^{+}$, and the quadratic functional $\Pi_\kappa^\pm$ on $X$ is defined as
\begin{equation}
\Pi_\kappa^\pm(v):=\pm\frac{1}{\kappa}R_d (v)+b_u(v,v)-\kappa R_p( v), \quad  v \in X.
\end{equation}

Consider the following minimizing problem:
\begin{equation}
y_\kappa^\pm=\arg \min_{v \in X}\Pi_\kappa^\pm(v),
\end{equation}
and it can be recognized that
\begin{equation}
y_\kappa^\pm=\frac{1}{2}\left(\kappa \epsilon_p \mp \frac{1}{\kappa} \epsilon_d\right),
\end{equation}
where $\epsilon_p,\epsilon_d \in X$ are the solutions of the following 'residual' problems:
\begin{equation}\label{res1}
b_u^S(\epsilon_p,v)=R_p(v),\quad b_u^S(\epsilon_d,v)=R_d(v),\quad \forall v\in V.
\end{equation}
Then it follows that
\begin{equation}
\pm (Q(u)-Q(u_h))=\Pi_\kappa^\pm(\kappa (u-u_h))\geq \Pi_\kappa^\pm(y_\kappa^\pm) = -\frac{1}{4}b_u^S\left(\kappa \epsilon_p \mp \frac{1}{\kappa}\epsilon_d,\kappa \epsilon_p \mp \frac{1}{\kappa}\epsilon_d\right),
\end{equation}
which is in a similar form with the front part of Eq. \eqref{ded1}.

If the equilibrium fields for the primal and adjoint 'residual' problems in Eq. \eqref{res1} are $\hat{\sigma}^{res}_p$ and $\hat{\sigma}^{res}_d$ that can be induced by a lower-order 'stress' bilinear form $b_{\sigma}^S(\cdot,\cdot)$, the bounding property of CRE gives
\begin{equation}\label{nsupb}
b_u^S\left(\kappa \epsilon_p \mp \frac{1}{\kappa}\epsilon_d,\kappa \epsilon_p \mp \frac{1}{\kappa}\epsilon_d\right)\leq b_{\sigma}^S\left(\kappa \hat{\sigma}^{res}_p \mp \frac{1}{\kappa}\hat{\sigma}^{res}_d,\kappa \hat{\sigma}^{res}_p \mp \frac{1}{\kappa}\hat{\sigma}^{res}_d\right),
\end{equation}
which is a natural result of the principle of minimum complementary energy (or called dual variational principle). Then similar bounds with those in Eq. \eqref{Ibounds} can be derived as
\begin{equation}
\left| Q(u)-Q(u_h)-\frac{1}{2}b_{\sigma}^S(\hat{\sigma}^{res}_p,\hat{\sigma}^{res}_d)\right|\leq \frac{1}{2}\sqrt{ b_{\sigma}^S(\hat{\sigma}^{res}_p,\hat{\sigma}^{res}_p)\cdot b_{\sigma}^S(\hat{\sigma}^{res}_d,\hat{\sigma}^{res}_d) },
\end{equation}
with $\kappa$ being taken as $\sqrt[4]{b_{\sigma}^S(\hat{\sigma}^{res}_d,\hat{\sigma}^{res}_d)/b_{\sigma}^S(\hat{\sigma}^{res}_p,\hat{\sigma}^{res}_p)}$.

Note that in the symmetric case in Section 2, one has
\begin{equation}
\begin{split}
& a_u=a^S_u\,,\;a_\sigma=a_\sigma^S\,,\\
& \hat{\sigma}^{res}_p=\hat{\sigma}-K(A(u_h))\,,\;\hat{\sigma}^{res}_d=\hat{\tilde{\sigma}}-K(A(\tilde{u}_h))\,,\\
& e_{CRE}(u_h,\hat{\sigma})=\|\hat{\sigma}^{res}_p\|_\sigma\,,\;e_{CRE}(\tilde{u}_h,\hat{\tilde{\sigma}})=\|\hat{\sigma}^{res}_d\|_\sigma\,.
\end{split}
\end{equation}
As discussed in Subsection 2.3, the bounding property of CRE was also identified as a consequence of minimizing complementary energy for the 'residual' problem. Therefore, in the sense of the principle of minimum complementary energy, the present bounding technique of goal-oriented error estimation for cases of non-symmetric bilinear forms can be considered as an extension of the CRE defined in symmetric cases.

\section{Goal-oriented error estimation for static response sensitivity analysis}

For the static response sensitivity analysis \cite{sensop} of linear structural systems, the first-order perturbation method is usually used to evaluate variations of response variables around their mean values resulting from the varying inputs. In the perturbed formulation of various variables, derivatives with respect to the input parameters are required, and those for the static responses are derived based on the finite element analysis at the central values through the perturbation method. Thus, the finite element descritization error propagates through the numerical results of the sensitivity derivatives of response variables with respect to the inputs, which will be evaluated by the constitutive relation error in this section.

\subsection{Primal problem for the first-order perturbation}

Suppose the description of the structural system is governed by several basic input parameters, one of which\footnote{Practically, a complex parameterized variational problem is involved due to the variation of a set of inputs. However, the sensitivity derivative with respect to each parameter can usually be considered independently. Thus only the first-order perturbation with respect to a single parameter is discussed in this section.} is denoted by $\beta$ with mean value $\bar{\beta}$. In this paper, only the input parameters describing the material properties and load variables are under consideration, i.e.
\begin{equation}
K=K(\beta)\,,\; f=f(\beta)\,,
\end{equation}
and the corresponding sensitivity to these input parameters is analyzed. The cases with basic geometrical parameters can be transformed into a similar form to those with material or load parameters, as stated in \textbf{Remark 1}.

Throughout the remainder of this paper, the following symbols are employed to represent the quantities for the first-order perturbation:
\begin{equation}
(\bullet):=(\bullet)(\bar{\beta}), \quad (\bullet)':=[\partial_{\beta}(\bullet)]|_{\beta=\bar{\beta}}.
\end{equation}
Moreover, the bilinear forms $a_{u}(\cdot,\cdot)$, $a_{{\sigma}}(\cdot,\cdot)$, the (semi-) norms $\|\cdot\|_{{u}}$, $\|\cdot\|_{{\sigma}}$ and the expression of $e_{CRE}$ represent the corresponding functionals when $K=K(\bar{\beta})$.
Besides, a bilinear form $a'_{{u}}(\cdot,\cdot): \mathcal{V}\times \mathcal{V}\mapsto \mathbb{R}$ associated with the derivatives with respect to the basic parameter $\beta$ is defined by
\begin{equation}
a'_{{u}}(u,v)=_{\mathcal{E}^*}\langle K'(Au),Av \rangle_{\mathcal{E}}\,,\;(u,v)\in \mathcal{V}\times\mathcal{V}\,.
\end{equation}

For notation, more spaces
\begin{equation}
\mathcal{X}=\mathcal{V}\times\mathcal{V}\,, \;\mathcal{Z}=\mathcal{E}^*\times\mathcal{E}^*\,,\;
\mathcal{X}_{h}=\mathcal{V}_{h}\times\mathcal{V}_{h}\,,
\end{equation}
are introduced. Then, two bilinear forms $\mathcal{A}_{u}(\cdot,\cdot): \mathcal{X}\times \mathcal{X} \to \mathbb{R}$ and
$\mathcal{A}_{{\sigma}}(\cdot,\cdot): \mathcal{Z }\times \mathcal{Z} \to \mathbb{R}$ are given as
\begin{equation}
\begin{split}
&  \mathcal{A}_{{u}}(\{{u},{U}\},\{{v},{V}\})=
a_{{u}}({u},{v})+a_{{u}}({U},{V})+\xi a'_{{u}}({u},{V})\,,\;(\{{u},{U}\},\{{v},{V}\})\in\mathcal{X}\times\mathcal{X}\,,\\
&  \mathcal{A}_{{\sigma}}(\{{\sigma},{\Sigma}\},\{{\tau},{\Gamma}\})=
a_{{\sigma}}({\sigma},{\tau})+a_{{\sigma}}({\Sigma},{\Gamma})+\xi a_{{\sigma}}(K'(K^{-1}{\sigma}),{\Gamma})\,,\;(\{{\sigma},{\Sigma}\},\{{\tau},{\Gamma}\})\in\mathcal{Z}\times\mathcal{Z}\,,
\end{split}
\end{equation}
where $\xi \in \mathbb{R}^{+}$ is a parameter to ensure that $\mathcal{A}_{{u}}$ is coercive and the quadratic functional $\mathcal{B}_{\sigma}(\cdot):\mathcal{Z}\to \mathbb{R},\{\tau,\Gamma\}\mapsto \mathcal{A}_{\sigma}(\{\tau,\Gamma\},\{\tau,\Gamma\})$ is positive definite (see \textbf{Remark 2}). It is obvious that $\mathcal{A}_{{u}}$ and $\mathcal{A}_{{\sigma}}$ are non-symmetric. The symmetric parts of $\mathcal{A}_u$ and $\mathcal{A}_\sigma$ are denoted by $\mathcal{A}_{{u}}^S$ and $\mathcal{A}_{{\sigma}}^S$, respectively.

The weak form of the primary problem at mean value of input parameter $\beta$ is given as: find ${u}\in \mathcal{V}$ such that
\begin{equation}\label{zero}
a_{{u}}({{u},{v}})=\langle f,v \rangle \quad \forall {v}\in \mathcal{V}\,.
\end{equation}
Differentiation of Eq.~\eqref{zero} with respect to $\beta$ gives the first-order perturbed equation as: find ${u}'\in \mathcal{V}$ such that
\begin{equation}\label{one}
a_{{u}}({{u}',{V}})=\langle f',V \rangle- a'_{{u}}({u},{V})\quad\forall {V}\in \mathcal{V}\,.
\end{equation}
Eqs. \eqref{zero} and \eqref{one} can be rewritten in a compact form as: find $\{{u},{U}\}\in \mathcal{X}$, where ${U}:=\xi{u}'$, such that
\begin{equation}\label{dprim}
\mathcal{A}_{{u}}(\{{u},{U}\},\{{v},{V}\})=\langle f,v \rangle + \xi\langle f',V \rangle \quad \forall \{{v},{V}\} \in \mathcal{X}\,.
\end{equation}


Adopting the finite element space $\mathcal{P}_{h}$, the finite element solution to this problem can be stated as: find $\{{u}_{h},{U}_{h}\}\in \mathcal{X}_{h}$ such that
\begin{equation}\label{FEdp}
\mathcal{A}_{{u}}(\{{u}_h,{U}_h\},\{{v},{V}\})=\langle f,v \rangle + \xi\langle f',V \rangle \quad \forall \{{v},{V}\} \in \mathcal{X}_h\,,
\end{equation}
with $U_h=\xi u'_h$. $\{{\sigma}_{h},{\Sigma}_{h}\}$, the finite element approximation of $\{{\sigma},{\Sigma}\}$, with $\Sigma:=\xi(\sigma'-K'(A(u)))$, is then obtained via the constitutive relation
\begin{equation}
   \{{\sigma}_h,{\Sigma}_h\}=\{K(A(u_h)),K(A(U_h))\}
\end{equation}
in the sense of distribution.

\vspace{4mm}
\noindent \textbf{Remark 1:} In this remark, a 3D problem in linear elasticity is taken as an example to show how to transform geometrical parameters to material-like parameters. Without loss of generality, the domain $\Omega$ that the structure occupies can be divided into several non-overlapped subdomains $\omega_{1},\cdots,\omega_{k},\cdots$, i.e. $\bigcup_{k}\omega_{k}=\Omega$, and a transformation $\mathcal{Y}_{k}$ can be defined for each subdomain $\omega_{k}$ to map it onto the standard domain $V=(-1,1)^{3}$, i.e. $\mathcal{Y}_{k}:\omega_{k}\to V,\boldsymbol{x}\mapsto\boldsymbol{\eta}$. Then the bilinear form $a_{{u}}$ can be represented by
\begin{equation}
\begin{split}
a_{{u}}(\boldsymbol{u},\boldsymbol{v}) & =\int_{\Omega} \nabla_{\boldsymbol{x}}\boldsymbol{u}:\boldsymbol{H}:\boldsymbol{v}\nabla_{\boldsymbol{x}}  \\
&=\sum_{\omega_{k}}\int_{V}(\nabla_{\boldsymbol{x}}\boldsymbol{\eta}\cdot\nabla_{\boldsymbol{\eta}}\boldsymbol{u}):\boldsymbol{H}:(\boldsymbol{u}\nabla_{\boldsymbol{\eta}}\cdot\boldsymbol{\eta}\nabla_{\boldsymbol{x}})\left|\det\nabla_{\boldsymbol{\eta}}\boldsymbol{x}\right| \\&=\sum_{\omega_{k}}\int_{V}\nabla_{\boldsymbol{\eta}}\boldsymbol{u}:\check{\boldsymbol{H}}:\boldsymbol{u}\nabla_{\boldsymbol{\eta}}\,,
\end{split}
\end{equation}
where $\boldsymbol{H}$ is Hooke's stiffness tensor, and $\check{\boldsymbol{H}}:=\boldsymbol{\eta}\nabla_{\boldsymbol{x}}\cdot\boldsymbol{H}\left|\det\nabla_{\boldsymbol{\eta}}\boldsymbol{x}\right|\cdot\nabla_{\boldsymbol{x}}\boldsymbol{\eta}$. It can be seen that the geometrical parameters for the structural system are all included in the 'equivalent' stiffness tensor  $\check{\boldsymbol{H}}$, so the cases with geometrical parameters can be treated as ones with material parameters. A similar treatment can be adopted for the loading functional $\langle f,\cdot\rangle$.

\vspace{2mm}
\noindent \textbf{Remark 2:} The determination of $\xi$ for 3D problems in linear elasticity is introduced in this remark. According to Hooke's Law for isotropic elastic material, the tensor $\boldsymbol{H}$ is represented in the index form as
\begin{equation}
H_{ijkl}=\lambda \delta_{ij}\delta_{kl}+\mu (\delta_{ik}\delta_{jl}+\delta_{il}\delta_{jk}),\;i,j,k,l=1,2,3,
\end{equation}
where $\lambda$ and $\mu$ are Lam\'{e} constants, satisfying $\mu >0$ and $3\lambda+2\mu>0$, and $\delta_{ij}$ is the Kronecker-delta (or unit) tensor. Moreover, one can obtain that
\begin{equation}
\int_{\Omega} \boldsymbol{\varepsilon}:\boldsymbol{H}:\boldsymbol{\varepsilon}=\lambda\|\mathrm{tr}\boldsymbol{\varepsilon}\|_{L^{2}}^{2}+2\mu\|\boldsymbol{\varepsilon}\|_{L^{2}}^{2}\geq \alpha \|\boldsymbol{\varepsilon}\|_{L_{2}}^{2}\geq 0,
\end{equation}
where $\alpha:=\min(2\mu,2\mu+3\lambda)>0$ and
\begin{equation}
\left(\sum_{i=1}^{3}\varepsilon_{ii}\right)^2\leq \left(\sum_{i=1}^{3}|\varepsilon_{ii}|\right)^2
\leq 3 \sum_{i=1}^{3}\varepsilon_{ii}^2\leq 3\boldsymbol{\varepsilon}:\boldsymbol{\varepsilon}
\end{equation}
is taken into consideration.
In addition, one has
\begin{equation}
\int_{\Omega}\boldsymbol{\varepsilon}:\boldsymbol{H}':\boldsymbol{\epsilon}  \geq -\frac{1}{2}\int_{\Omega}\sum_{i,j,k,l=1}^{3} H'_{ijkl}\left(\varepsilon_{ij}^{2}+\epsilon_{kl}^{2}\right)
\geq -\frac{1}{2} \psi \left(\|\boldsymbol{\varepsilon}\|_{L^{2}}^{2}+\|\boldsymbol{\epsilon}\|_{L^{2}}^{2}\right),
\end{equation}
where $\psi:=\max_{i,j,k,l=1,2,3}|H'_{ijkl}|$. Therefore, when $0<\xi<2\alpha/\psi$,
\begin{equation}
\int_{\Omega}\boldsymbol{\varepsilon}:\boldsymbol{H}:\boldsymbol{\varepsilon} +\int_{\Omega}\boldsymbol{\epsilon}:\boldsymbol{H}:\boldsymbol{\epsilon} +
\int_{\Omega}\boldsymbol{\varepsilon}:\boldsymbol{H}':\boldsymbol{\epsilon}
\geq \left(\alpha-\frac{1}{2}\xi \psi\right)\left( \|\boldsymbol{\varepsilon}\|_{L^{2}}^{2}+\|\boldsymbol{\epsilon}\|_{L^{2}}^{2} \right)\geq 0\,.
\end{equation}

Then one can conclude from Korn's inequality (see \cite{introfunc}) that
\begin{equation}
\begin{split}
\mathcal{A}_{u}(\{\boldsymbol{v},\boldsymbol{V}\},\{\boldsymbol{v},\boldsymbol{V}\}) & =\int_{\Omega}\boldsymbol{\varepsilon}(\boldsymbol{v}):\boldsymbol{H}:\boldsymbol{\varepsilon}(\boldsymbol{v}) +\int_{\Omega}\boldsymbol{\varepsilon}(\boldsymbol{V}):\boldsymbol{H}:\boldsymbol{\varepsilon}(\boldsymbol{V}) +
\int_{\Omega}\boldsymbol{\varepsilon}(\boldsymbol{v}):\boldsymbol{H}':\boldsymbol{\varepsilon}(\boldsymbol{V})\\
& \geq \left(\alpha-\frac{1}{2}\xi \psi\right)\left( \|\boldsymbol{\varepsilon}(\boldsymbol{v})\|_{L^{2}}^{2}+\|\boldsymbol{\varepsilon}(\boldsymbol{V})\|_{L^{2}}^{2} \right)\\
& \geq C_1\left(\alpha-\frac{1}{2}\xi \psi\right)\left( \|\boldsymbol{v}\|_{\mathcal{V}}^2+\|\boldsymbol{V}\|_{\mathcal{V}}^{2} \right)\quad \forall \{\boldsymbol{v},\boldsymbol{V}\}\in \mathcal{X}\,,
\end{split}
\end{equation}
where $\boldsymbol{\varepsilon}(\bullet)=(\nabla\bullet+\bullet\nabla)/2$, $C_1$ is a positive constant, $\mathcal{V}=\{\boldsymbol{v}\in [H^1(\Omega)]^3:\boldsymbol{v}|_{\partial\Omega_{D}}=\boldsymbol{0}\}$, $\Omega\in\mathbb{R}^3$ is the problem domain with Lipschitz boundary $\partial\Omega$, $\partial\Omega_{D}\subset \partial\Omega$ is the Dirichlet boundary, and $\partial\Omega_{D}\neq\emptyset$. Thus $\mathcal{A}_u$ is coercive. It follows that
\begin{equation}
\begin{split}
& \mathcal{A}_{\sigma}(\{\boldsymbol{\tau},\boldsymbol{\Gamma}\},\{\boldsymbol{\tau},\boldsymbol{\Gamma}\})\\
=& \int_{\Omega}(\boldsymbol{H}^{-1}:\boldsymbol{\tau}):\boldsymbol{H}:(\boldsymbol{H}^{-1}:\boldsymbol{\tau}) +\int_{\Omega}(\boldsymbol{H}^{-1}:\boldsymbol{\Gamma}):\boldsymbol{H}:(\boldsymbol{H}^{-1}:\boldsymbol{\Gamma}) +
\int_{\Omega}(\boldsymbol{H}^{-1}:\boldsymbol{\Gamma}):\boldsymbol{H}':(\boldsymbol{H}^{-1}:\boldsymbol{\Gamma})\\
\geq & \left(\alpha-\frac{1}{2}\xi \psi\right)\left( \|\boldsymbol{H}^{-1}:\boldsymbol{\tau}\|_{L^{2}}^{2}+\|\boldsymbol{H}^{-1}:\boldsymbol{\Gamma}\|_{L^{2}}^{2} \right)>0
\quad \forall \{\boldsymbol{\tau},\boldsymbol{\Gamma}\}\in \mathcal{Z}\setminus \{\boldsymbol{0},\boldsymbol{0}\}\,,
\end{split}
\end{equation}
i.e. $\mathcal{B}_{\sigma}(\cdot):\mathcal{Z}\to \mathbb{R},\{\boldsymbol{\tau},\boldsymbol{\Gamma}\}\mapsto \mathcal{A}_{\sigma}(\{\boldsymbol{\tau},\boldsymbol{\Gamma}\},\{\boldsymbol{\tau},\boldsymbol{\Gamma}\})$ is positive definite.

\subsection{Error estimator extended from CRE}

In this case, bilinear forms $\mathcal{A}_u$ and $\mathcal{A}_{\sigma}$ are non-symmetric, which can be treated using the technique proposed in Section \ref{nonsym}. As preciously introduced, a residual linear functional $\mathcal{R}_p \in \mathcal{X}^*$ is defined as
\begin{equation}
\mathcal{R}_p(\{{v},{V}\})=
\langle f,v \rangle + \xi\langle f',V \rangle-\mathcal{A}_{{u}}(\{{u}_h,{U}_h\},\{{v},{V}\})\,,\; \{{v},{V}\} \in \mathcal{X}\,.
\end{equation}
A statically admissible field pair for the residual problem, $\{\hat{{\sigma}}^{res}, \hat{{\Sigma}}^{res}\}\in \mathcal{Z}$ that satisfies
\begin{equation}\label{dperes}
\mathcal{A}_{{\sigma}}^S(\{\hat{{\sigma}}^{res}, \hat{{\Sigma}}^{res}\},\{K(A(v)),K(A(V))\})=\mathcal{R}_p(\{{v},{V}\})\quad \forall \{{v},{V}\} \in \mathcal{X}\,,
\end{equation}
can be further obtained. Then an error estimator is defined by the admissible field pair as
\begin{equation}
{E}(\hat{{\sigma}}^{res},\hat{{\Sigma}}^{res})
:= \sqrt{\mathcal{A}_{{\sigma}}^{S}(\{\hat{{\sigma}}^{res},\hat{{\Sigma}}^{res}\},\{\hat{{\sigma}}^{res},\hat{{\Sigma}}^{res}\})}\,,
\end{equation}
which can be considered as an extension of CRE, as stated in Section 3. This estimator has the bounding property given in Eq. \eqref{nsupb}.



\subsection{Goal-oriented error estimation associated with the sensitivity derivative fields of displacements}

Defined via a linear bounded functional $J \in \mathcal{V}^*$ associated with the derivative solution field ${u}'$, the quantity of interest is denoted by $J(u')=\frac{1}{\xi}J(U)$ and the computed value of the quantity is represented as $J(u'_h)=\frac{1}{\xi}J(U_h)$. Then an adjoint problem can be defined for the output $J(u')$, and its corresponding weak form is given as: find $\{{{w}},{{W}}\}\in \mathcal{X}$ such that
\begin{equation}\label{dadj}
\mathcal{A}_{{u}}(\{{v},{V}\},\{{{w}},{{W}}\})=\frac{1}{\xi}J({V})\quad\forall \{{v},{V}\} \in \mathcal{X}\,,
\end{equation}
and it can be discretized by finite elements as: find $\{{{w}}_{h},{{W}}_{h}\}\in \mathcal{X}_{h}$ such that
\begin{equation}\label{FEdd}
\mathcal{A}_{{u}}(\{{v},{V}\},\{w_h,W_h\})=\frac{1}{\xi}J({V})\quad\forall \{{v},{V}\} \in \mathcal{X}_h\,.
\end{equation}

Analogous to Subsection 2.4, the computed error in quantity $J$ can be represented as
\begin{equation}\label{derrorrep}
 J(u')-J(u'_h)=\mathcal{A}_{{u}}(\{{u}-{u}_{h},{U}-{U}_{h}\},\{{{w}}-{{w}}_{h},{{W}}-{{W}}_{h}\}).
\end{equation}
After defining the residual linear functional $\mathcal{R}_d\in \mathcal{X}^*$ for the adjoint problem as
\begin{equation}
\mathcal{R}_d(\{{v},{V}\})=
\frac{1}{\xi}J(V)-\mathcal{A}_{{u}}(\{{v},{V}\},\{{w}_h,{W}_h\})\,,\; \{{v},{V}\} \in \mathcal{X}\,,
\end{equation}
the statically admissible field pair $\{\hat{{\tau}}^{res}, \hat{{\Gamma}}^{res}\}\in \mathcal{Z}$ for the residual problem of adjoint problem satisfying
\begin{equation}\label{dderes}
\mathcal{A}_{{\sigma}}^S(\{\hat{{\tau}}^{res}, \hat{{\Gamma}}^{res}\},\{K(A(v)),K(A(V))\})=\mathcal{R}_d(\{{v},{V}\})\quad \forall \{{v},{V}\} \in \mathcal{X}\,,
\end{equation}
is then obtained.

As presented in Section 3, strict upper and lower bounds are given by the estimator extended from CRE as
\begin{equation}\label{dbounds}
\left|J(u')-J(u'_h)-\frac{1}{2}\mathcal{A}^S_{{\sigma}} (\{\hat{{\sigma}}^{res},\hat{{\Sigma}}^{res}\} ,\{\hat{{{\tau}}}^{res},\hat{{{\Gamma}}}^{res}\} )    \right|
 \leq  \frac{1}{2}{E}(\hat{{\sigma}}^{res},\hat{{\Sigma}}^{res})
 \cdot{E}(\hat{{{\tau}}}^{res},\hat{{{\Gamma}}}^{res})\,.
\end{equation}
Then the corresponding strict upper and lower bounds for quantity $J(u')$ are, respectively,
\begin{equation}
\begin{split}
& J^{upper}=J(u'_h)+\frac{1}{2}\mathcal{A}^S_{{\sigma}} (\{\hat{{\sigma}}^{res},\hat{{\Sigma}}^{res}\} ,\{\hat{{{\tau}}}^{res},\hat{{{\Gamma}}}^{res}\} )+\frac{1}{2}{E}(\hat{{\sigma}}^{res},\hat{{\Sigma}}^{res})
 \cdot{E}(\hat{{{\tau}}}^{res},\hat{{{\Gamma}}}^{res})\,,\\
& J^{lower}=J(u'_h)+\frac{1}{2}\mathcal{A}^S_{{\sigma}} (\{\hat{{\sigma}}^{res},\hat{{\Sigma}}^{res}\} ,\{\hat{{{\tau}}}^{res},\hat{{{\Gamma}}}^{res}\} )-\frac{1}{2}{E}(\hat{{\sigma}}^{res},\hat{{\Sigma}}^{res})
 \cdot{E}(\hat{{{\tau}}}^{res},\hat{{{\Gamma}}}^{res})\,.
\end{split}
\end{equation}

\subsection{Quantities of response sensitivity derivatives}

Let $g: \mathcal{V}\to \mathbb{R}$ denote a linear or non-linear functional with respect to the displacement field solution ${u}={u}(\beta)$ with the basic parameter $\beta$, representing a scalar-valued quantity in the static response of a structural system. This response quantity $G=g({u})$ can also be considered as a function with respect to $\beta$
\begin{equation}
G=g({u}(\beta)).
\end{equation}

Usually, the quantity of interest is the sensitivity derivative of $G$, which can be expressed in the following form on the basis of the chain rule of derivatives:
\begin{equation}
G'=Dg[{u}]({u}')\,,
\end{equation}
where $Dg[{u}](\cdot)\in \mathcal{V}^*$ is the G\^{a}teaux derivative of the functional $g$, i.e.
\begin{equation}
Dg[{u}]({v}):=\lim_{t\rightarrow 0}\frac{g({u}+t{v})-g({u})}{t}\,,\; {v} \in \mathcal{V}\,.
\end{equation}

Since the G\^{a}teaux derivative is a linear bounded functional in $\mathcal{V}$, the quantity of interest $G'$ manifests itself in the global functional form  associated with the derivative field $u'$, i.e.
\begin{equation}
J(\bullet):=Dg[{u}](\bullet)\,,
\end{equation}
so the descritization error in $G'=J(u')$ can be estimated by the techniques introduced in this section.

\vspace{2mm}
\noindent \textbf{Remark 3:} The finite element approximation \eqref{FEdp} of first-order perturbation can be written in the matrix form as
\begin{equation}
\mathbf{Ku}=\mathbf{f}\,,\quad \mathbf{KU}=\xi\mathbf{f}'-\xi\mathbf{K}'\mathbf{u}\,,
\end{equation}
where $\mathbf{U}=\xi\mathbf{u}'$, $\mathbf{K}$, $\mathbf{f}$ and $\mathbf{u}$ are the commonly defined global matrix of stiffness, global vectors of load and displacement, respectively, and $\mathbf{K}'$, $\mathbf{f}'$ and $\mathbf{u}'$ are their derivatives with respect to the parameter $\beta$. The computed value of the quantity of interest $J(u')$ can also be represented as
\begin{equation}
J(u'_h)=\frac{1}{\xi}\mathbf{g}^T \mathbf{U}\,,
\end{equation}
where $\mathbf{g}$ is referred to as an extracting vector, with its $i$th component $\mathbf{g}_i=J(N_i)$, $N_i$ being the shape function associated with the $i$th degree of freedom. This technique is also termed the direct differentiation method for sensitivity.

Alternatively, the adjoint state method, or referred to as the adjoint variable method \cite{sensop}, is based on an adjoint problem
\begin{equation}
\mathbf{K}^T\boldsymbol{\lambda}=\mathbf{g}\,,
\end{equation}
which is actually included in Eq. \eqref{FEdd}, and the output is computed as
\begin{equation}
J(u'_h)=\boldsymbol{\lambda}^T(\mathbf{f}'-\mathbf{K}'\mathbf{u})\,.
\end{equation}
It can be noted that
\begin{equation}
\boldsymbol{\lambda}^T(\mathbf{f}'-\mathbf{K}'\mathbf{u})=\mathbf{g}^T\mathbf{K}^{-1}(\mathbf{f}'-\mathbf{K}'\mathbf{u})=\mathbf{g}^T\mathbf{u}'\,,
\end{equation}
meaning that the same computed value for the output of interest can be obtained via either the direct differentiation method or the adjoint state method. Therefore, the proposed error estimation technique is also applicable to the latter for sensitivity analysis.

\section{Numerical examples for model problems}

\subsection{Beam model problem}
The proposed technique is exemplified by the Bernoulli-Euler beam model in this subsection. The transverse deflection $u$ of a beam is taken as the displacement field, curvature $\phi$ as the strain field, bending moment $M$ as the stress field. Distributed transverse load $q$ is applied on the beam. In the example, the operators $A=A^*=\partial_{xx}$ and $K=EI(x)\geq \min_{s\in Y}EI(s)>0$, where $x\in Y$ and $Y$ is the one-dimensional interval that the beam occupies. The space of admissible displacements $\mathcal{V}$ is a subspace of $H^2(Y)$ whose elements satisfy the homogeneous Dirichlet boundary conditions, and the space of strain or stress is $\mathcal{E}^*=\mathcal{E}=L^2(Y)$. The three types of equations governing the problem are listed as follows:
\begin{itemize}
  \item The compatibility conditions: $u\in \mathcal{V}$, including the conditions $u=0$ and/or $\partial_x u=0$ at the Dirichlet boundaries;
  \item The equilibrium conditions: $M\in \mathcal{E}^*$, $\int_X M\partial_{xx}v=\int_X qv$, including the conditions $M=0$ and/or $\partial_x M=0$ at the Neumann boundaries;
  \item The constitutive relation: $M=EI\partial_{xx}u$\,.
\end{itemize}

As stated in Section 4, if there exists a dimensionless parameter $\beta$ with its mean value $\bar{\beta}$, an approximate solution $\{u_h,U_h\}$ can be obtained for the field pair $\{u,\xi u'\}$ by finite element analysis \eqref{FEdp} under mesh size $h$. Similarly, finite element solution $\{w_h,W_h\}$ for the adjoint problem in \eqref{FEdd} can be obtained. Analogous to \textbf{Remark 2}, the parameter $\xi$ should satisfy $0<\xi<2/\max_{x\in Y}(EI'(x)/EI(x))$.

For the primal problem \eqref{dprim}, the 'equilibrated residual' field pair $\{\hat{M}^{res}_h,\hat{\mu}^{res}_h\}$ satisfying Eq.~\eqref{dperes} should be constructed, which can be explicitly rewritten as
\begin{equation}
\begin{split}
\int_Y\left(\hat{M}^{res}_h+\frac{\xi EI'}{2 EI}\hat{\mu}^{res}_h\right)\partial_{xx}v=\int_Y qv -\int_Y EI~\partial_{xx}u_h~ \partial_{xx}v \quad \forall v\in \mathcal{V}\,,&\\
\int_Y\left(\hat{\mu}^{res}_h+\frac{\xi EI'}{2 EI}\hat{M}^{res}_h\right)\partial_{xx}V=\xi\int_Y q'V -\int_Y EI~\partial_{xx}U_h~ \partial_{xx}V \quad & \\
-\xi \int_Y EI'~\partial_{xx}u_h~\partial_{xx}V \quad\forall V\in \mathcal{V}\,.&
\end{split}
\end{equation}
Let
\begin{equation}\label{trans}
\begin{split}
& \hat{M}_h=\hat{M}^{res}_h+\frac{\xi EI'}{2 EI}\hat{\mu}^{res}_h+EI~\partial_{xx}u_h\,,\\
& \hat{\mu}_h=\hat{\mu}^{res}_h+\frac{\xi EI'}{2 EI}\hat{M}^{res}_h+EI~\partial_{xx}U_h+\xi EI'~\partial_{xx}u_h\,.
\end{split}
\end{equation}
Then one can recognize that $\{\hat{M}_h,\hat{\mu}_h\}$ are statically admissible solutions to the following problems:
\begin{equation}
\begin{split}
& \int_Y\hat{M}_h~\partial_{xx}v=\int_Y qv  \quad \forall v\in \mathcal{V}\,,\\
& \int_Y\hat{\mu}_h~\partial_{xx}V=\xi\int_Y q'V  \quad \forall V\in \mathcal{V}\,,
\end{split}
\end{equation}
and they can be constructed based on $\{u_h,U_h\}$ using the recovery technique proposed in \cite{Guo1}. Employing the inverse transformation of \eqref{trans}, $\{\hat{M}^{res}_h,\hat{\mu}^{res}_h\}$ are obtained as
\begin{equation}
\begin{split}
& \hat{M}_h^{res}=\frac{1}{1-\lambda^2}\left[\left(\hat{M}_h-EI~\partial_{xx}u_h\right)-\lambda\left(\hat{\mu}_h-EI~\partial_{xx}U_h-\xi EI'~\partial_{xx}u_h\right)\right]\,,\\
& \hat{\mu}_h^{res}=\frac{1}{1-\lambda^2}\left[\left(\hat{\mu}_h-EI~\partial_{xx}U_h-\xi EI'~\partial_{xx}u_h\right)-\lambda\left(\hat{M}_h-EI~\partial_{xx}u_h\right)\right]\,,
\end{split}
\end{equation}
where $\lambda(x):=(\xi EI'(x))/(2EI(x))$. For the adjoint problem \eqref{dadj} with respect to quantity $J$, analogously, the 'equilibrated residual' field pair $\{\hat{T}^{res}_h,\hat{\psi}^{res}_h\}$ can be constructed from $\{w_h,W_h\}$.

Finally, the terms in \eqref{dbounds} are explicitly written as
\begin{equation}
\begin{split}
E^2(\hat{M}^{res}_h,\hat{\mu}^{res}_h) =&\int_Y\frac{(\hat{M}^{res}_h)^2}{EI}+\int_Y\frac{(\hat{\mu}^{res}_h)^2}{EI}+\xi\int_Y\frac{EI'}{EI}\frac{\hat{M}^{res}_h\hat{\mu}^{res}_h}{EI}\,,\\
E^2(\hat{T}^{res}_h,\hat{\psi}^{res}_h)
=&\int_Y\frac{(\hat{T}^{res}_h)^2}{EI}+\int_Y\frac{(\hat{\psi}^{res}_h)^2}{EI}+\xi\int_Y\frac{EI'}{EI}\frac{\hat{T}^{res}_h\hat{\psi}^{res}_h}{EI}\,,\\
\mathcal{A}_{\sigma}^S(\{\hat{M}^{res}_h,\hat{\mu}^{res}_h\},\{\hat{T}^{res}_h,\hat{\psi}^{res}_h\})
=&\int_Y\frac{\hat{M}^{res}_h \hat{T}^{res}_h}{EI}+\int_Y\frac{\hat{\mu}^{res}_h \hat{\psi}^{res}_h}{EI}\\
& +\frac{\xi}{2}\int_Y\frac{EI'}{EI}\frac{\hat{M}^{res}_h\hat{\psi}^{res}_h+\hat{T}^{res}_h\hat{\mu}^{res}_h}{EI}\,.
\end{split}
\end{equation}

\vspace{2mm}
\noindent \textbf{Remark 4:} In the beam model problem, which is a typical $C^1$ problem, the representation of quantity error $J(u')-J(u'_h)$ in Eq. \eqref{derrorrep} gives
\begin{equation}
|J(u')-J(u'_h)|\leq C_2 h^{2(p+1-2)}=C_2 h^{2p-2}\,,
\end{equation}
where $C_2$ is a positive constant independent of $h$, and $p$ is the interpolation order of the finite elements used.

The asymptotic behavior of the CRE, as shown in \cite{error1}, indicates that there exists a constant $C_3\geq 1$ independent of $h$ such that
\begin{equation}
e_{CRE}(u_h,\hat{\sigma}_h)\leq C_3~\|u-u_h\|_u\,,
\end{equation}
where $\hat{\sigma}_h$ is constructed from $u_h$ by the element equilibrium technique (EET, see \cite{error1}), a recovery technique. It is concluded that the quantity error $Q(u)-Q(u_h)$ and the bounding gap $Q^{upper}-Q^{lower}=e_{CRE}({u}_{h},\hat{{\sigma}}_{h})\cdot e_{CRE}(\tilde{{u}}_{h},\hat{\tilde{{\sigma}}}_{h})$ are of the same convergence rate as $h\rightarrow 0$. As an extension, the same convergence property can be achieved for the proposed error bounds of quantities with respect to sensitivity derivative fields.

\subsection{Numerical example 1: a portal frame}

\begin{figure}
  \begin{minipage}[t]{0.48\linewidth}
  \centering
   \includegraphics[width=.85\textwidth]{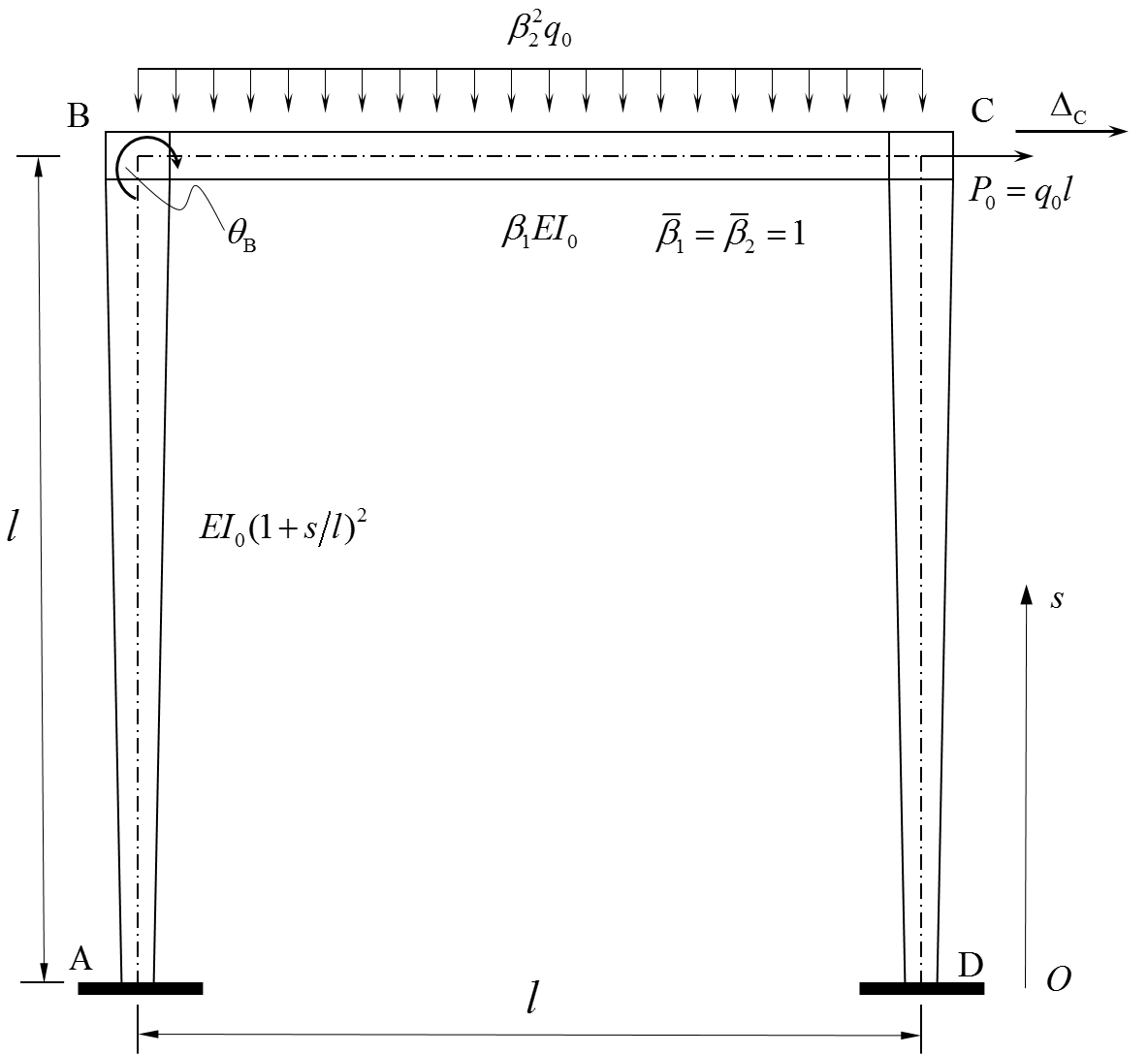}
  \caption{A portal frame with variable cross-sections of columns}\label{frame}
  \end{minipage}
  \hfill
  \begin{minipage}[t]{0.48\linewidth}
  \centering
   \includegraphics[width=\textwidth]{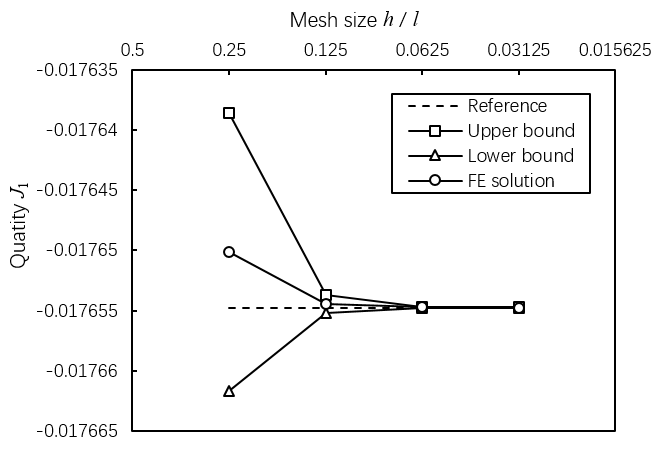}
  \caption{FE solutions, upper and lower bounds of $J_1$ ($\xi=1.0$) for the portal frame}\label{case1quan}
  \end{minipage}
\end{figure}

As shown in Figure \ref{frame}, a portal frame is under consideration. The flexible stiffness of each column, AB or DC, varies along the axis as $EI(s)=EI_0(1+s/l)^2$, where $l$ is the length of the column, and the flexion stiffness of the beam BC,  $\beta_1 EI_0$, is a constant. A uniformly distributed transverse load $\beta_2^2 q_0$ is applied on the beam and a horizontal concentrated load $P_0=q_0 l$ is prescribed at point C. $\beta_1$ and $\beta_2$ are two non-dimensional parameters, and $\bar{\beta}_1=\bar{\beta}_2=1$. The same uniform mesh made of third-order Hermitian beam elements $(p=3)$ is used in finite element analysis, and the size of a mesh is denoted by $h$, the length of an element. In the numerical examples, relative error (RE) of the finite element solution $J_h$ of a quantity $J$ is calculated as $RE(J_h)=|J-J_h|/|J|$, and that of the bounding gap is calculated as $RE(J^{upper},J^{lower})=(J^{upper}-J^{lower})/|J|$.

\textbf{Case 1:} $J_1=[\partial_{\beta_1}\underline{\Delta}_{\mathrm{C}}]_{(\bar{\beta}_1,\bar{\beta}_2)}$

$\underline{\Delta}_{\mathrm{C}}$ is the non-dimensionalized horizontal displacement at point C, i.e. $\underline{\Delta}_{\mathrm{C}}=\frac{\Delta_{\mathrm{C}}EI_0}{q_0 l^4}$. Under a refined mesh $h/l=1/50$, one has $\underline{\Delta}_{\mathrm{C}}(\bar{\beta}_1,\bar{\beta}_2)=0.0430866$, and the reference value of $J_1$ is found to be $-0.0176547$.

\begin{figure}
  \begin{minipage}[t]{0.48\linewidth}
  \centering
   \includegraphics[width=\textwidth]{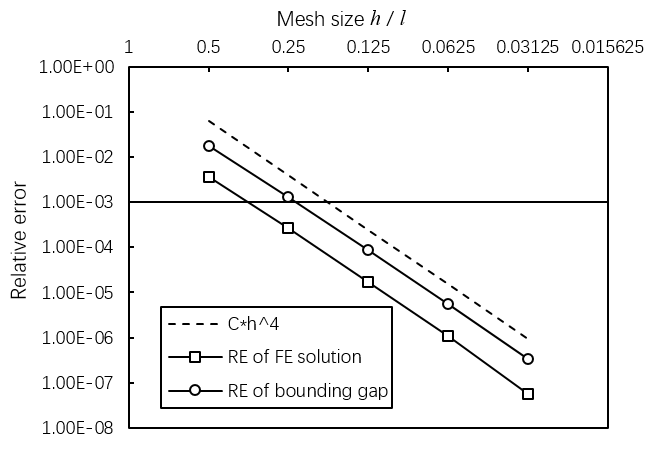}
  \caption{Relative errors of FE solutions and bounding gaps of $J_1$ ($\xi=1.0$) for the portal frame}\label{case1RE}
  \end{minipage}
  \hfill
  \begin{minipage}[t]{0.48\linewidth}
  \centering
   \includegraphics[width=\textwidth]{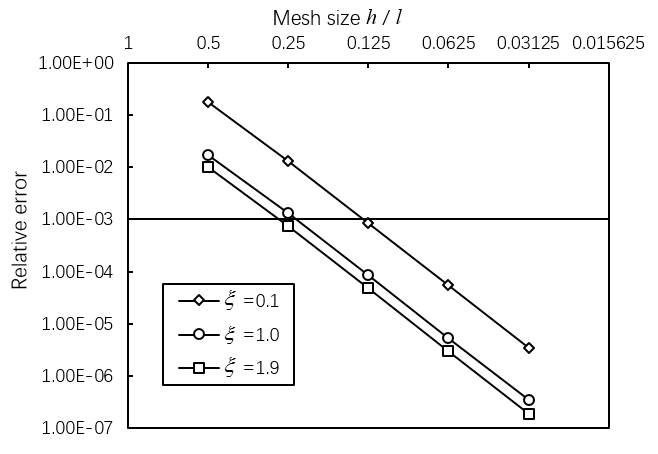}
  \caption{Relative errors of bounding gaps of $J_1$ for the portal frame with $\xi=0.1,~1.0,~1.9$}\label{case1RExi}
  \end{minipage}
\end{figure}

In this case, the range of $\xi$ is $0<\xi<2$. Taking $\xi=1.0$, it is seen in Figure \ref{case1quan} that strict bounding property of $J_1$ is achieved for various mesh densities. The relative errors of FE solutions and bounding gaps of the quantity are shown in Figure \ref{case1RE}. The relative error of bounding gap will be less than 0.1\% when the mesh size $h/l$ is smaller than $0.25$, and both the finite element solution and the bounding gap have the same convergence rate of $h^4$, showing super-convergent asymptotic property of the proposed error bounds.

When taking different values of $\xi$, the relative errors of bounding gaps versus decreased mesh size with $\xi=0.1,~1.0,~1.9$ are plotted in Figure \ref{case1RExi}. It is observed that the convergence rate of bounding gap keeps unchanged when the value of $\xi$ varies in the admissible range.

\textbf{Case 2:} $J_2=[\partial_{\beta_2}\underline{\theta}_{\mathrm{B}}]_{(\bar{\beta}_1,\bar{\beta}_2)}$

$\underline{\theta}_{\mathrm{B}}$ is the non-dimensionalized slope at point B, i.e. $\underline{\theta}_{\mathrm{B}}=\frac{\theta_{\mathrm{B}}EI_0}{q_0 l^3}$. Under a refined mesh $h/l=1/50$, one has $\underline{\theta}_{\mathrm{B}}(\bar{\beta}_1,\bar{\beta}_2)=0.0444687$, and the reference value of $J_2$ is found to be $0.0122243$.

\begin{figure}
  \begin{minipage}[t]{0.48\linewidth}
  \centering
   \includegraphics[width=\textwidth]{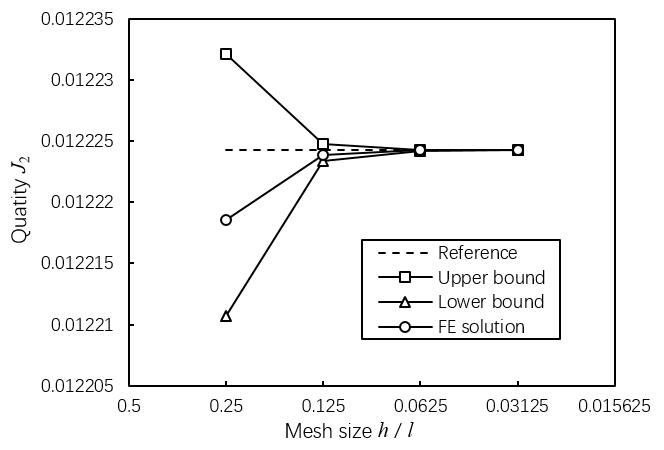}
  \caption{FE solutions, upper and lower bounds of $J_2$ for the portal frame}\label{case2quan}
  \end{minipage}
  \hfill
  \begin{minipage}[t]{0.48\linewidth}
  \centering
   \includegraphics[width=\textwidth]{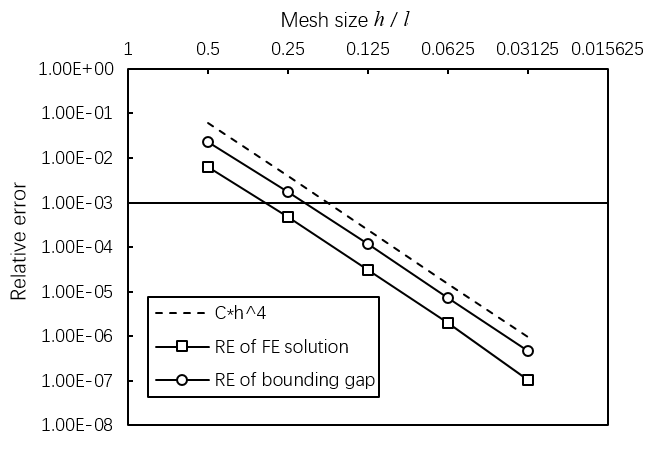}
  \caption{Relative errors of FE solutions and bounding gaps of $J_2$ for the portal frame}\label{case2RE}
  \end{minipage}
\end{figure}

Since the load is exclusively parameterized, $\xi$ is irrelevant in this case. The numerical results for quantity $J_2$ are illustrated in Figure \ref{case2quan}, assessing the strict bounding property of the proposed
goal-oriented error estimation for sensitivity derivative again. Figure \ref{case2RE} shows the relative errors of FE solutions and bounding gaps of the quantity $J_2$, both having the convergence rate of $h^4$. Super-convergence has been achieved by the proposed bounding gap.

\subsection{Numerical example 2: a membrane on an elastic foundation}

\begin{figure}[hbt]
  \begin{minipage}[t]{0.48\linewidth}
  \centering
   \includegraphics[width=0.85\textwidth]{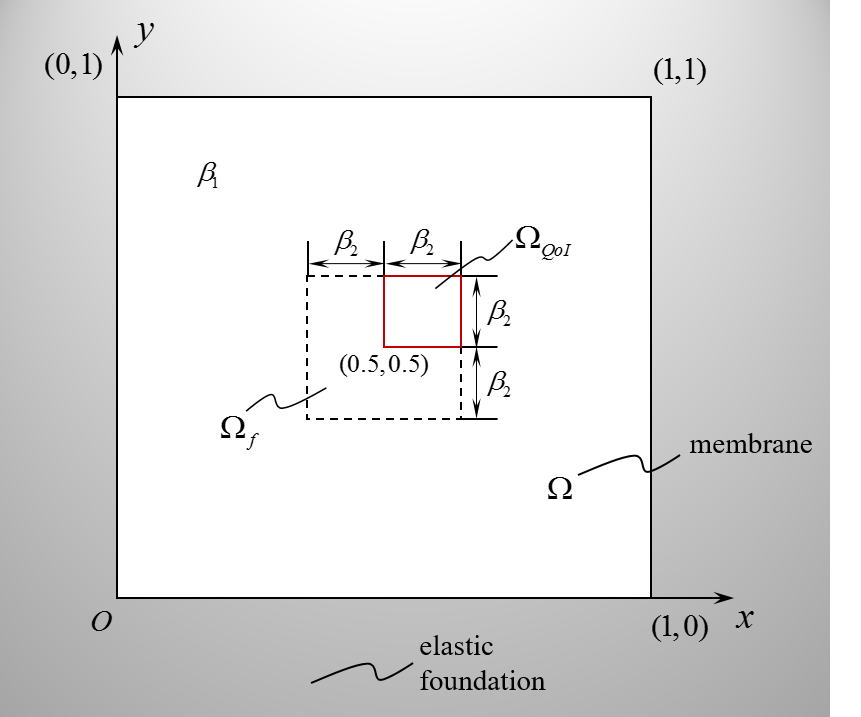}
  \caption{A membrane on an elastic foundation}\label{mem}
  \end{minipage}
  \hfill
  \begin{minipage}[t]{0.48\linewidth}
  \centering
   \includegraphics[width=\textwidth]{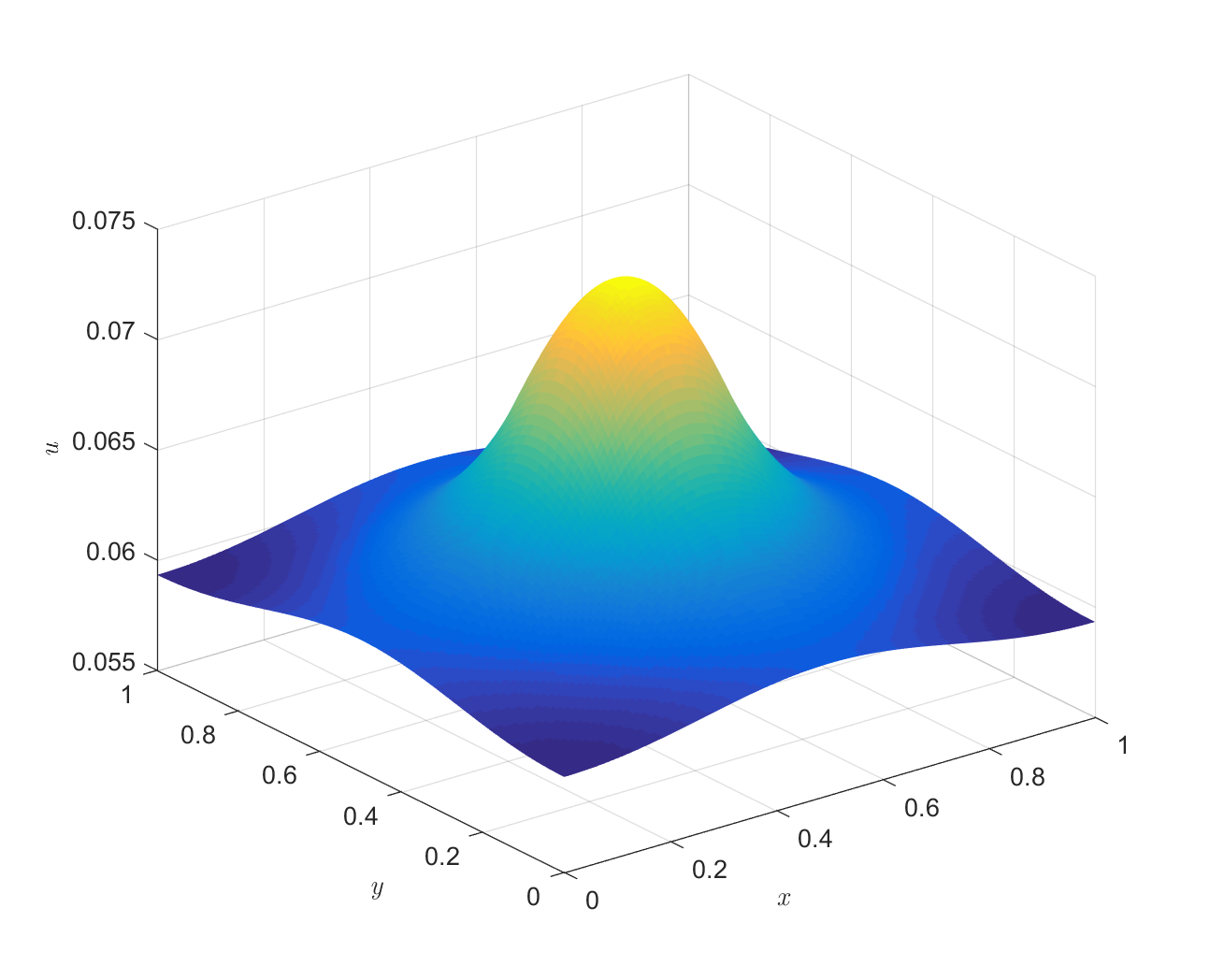}
  \caption{Finite element approximation of $u$ under a refined mesh $h=1/128$ in the membrane problem}\label{disp}
  \end{minipage}
\end{figure}

A membrane on an elastic foundation is considered in this subsection, which can be classified into the abstract formulation in Subsection 2.1 as well. As shown in Figure \ref{mem}, a square membrane, defined in the domain $\Omega=(0,1)^2\subset \mathbb{R}^2$ with the free boundary $\partial \Omega$, is settled on an elastic foundation with linear reaction with respect to the deflection. The governing equation of this problem is given in the weak form as: find $u\in \mathcal{V}=H^1(\Omega)$ such that
\begin{equation}
\int_\Omega a\nabla u\cdot\nabla v +\int_\Omega kuv =\int_\Omega fv\quad \forall v \in \mathcal{V}\,,
\end{equation}
where $a(\boldsymbol{x})=\beta_1$ and $k(\boldsymbol{x})=1$, $\boldsymbol{x}\in\Omega$, denote the stiffness of the membrane and the foundation, respectively, and the distributed load $f$ is given as
\begin{equation*}
f(\boldsymbol{x})=
\left\{
  \begin{array}{ll}
    1, & \boldsymbol{x}\in \Omega_f=(0.5-\beta_2,0.5+\beta_2)^2\,, \\
    0, & \mathrm{otherwise}\,,
  \end{array}
\right.
\quad \boldsymbol{x}\in \Omega\,.
\end{equation*}
The mean values of the two input parameters $\beta_1$ and $\beta_2$ are set to be $\bar{\beta}_1=1$ and $\bar{\beta}_2=1/8$, respectively.

Uniform meshes with bilinear quadrilateral elements are adopted in the finite element analysis of this problem, and the characterized mesh size, i.e. the length of a side of an element, is denoted by $h$. The finite element approximation of $u$ at the mean values $(\bar{\beta}_1,\bar{\beta}_2)$ under a refined mesh $h=1/128$ is shown in Figure \ref{disp}.

In this example, the derivatives of average displacement in the domain $\Omega_{QoI}=(0.5,0.625)^2$ with respect to the two parameters are considered as the quantities of interest, i.e. $J_i=\partial_{\beta_i}\left[\left(\int_{\Omega_{QoI}}u\right)/\left|\Omega_{QoI}\right|\right]$\footnote{For simplification, the fact that quantities $J_i$, $i=1,2$, are computed at the mean values of the parameters is not written explicitly in this example.}, $i=1,2$, $\left|\Omega_{QoI}\right|=1/64$ being the area of $\Omega_{QoI}$.

\textbf{Case 1:} $J_1=\partial_{\beta_1}\left[\left(\int_{\Omega_{QoI}}u\right)/\left|\Omega_{QoI}\right|\right]=64\int_{\Omega_{QoI}}\partial_{\beta_1}u$

\begin{figure}[htb]
  \centering
  \includegraphics[scale=0.3]{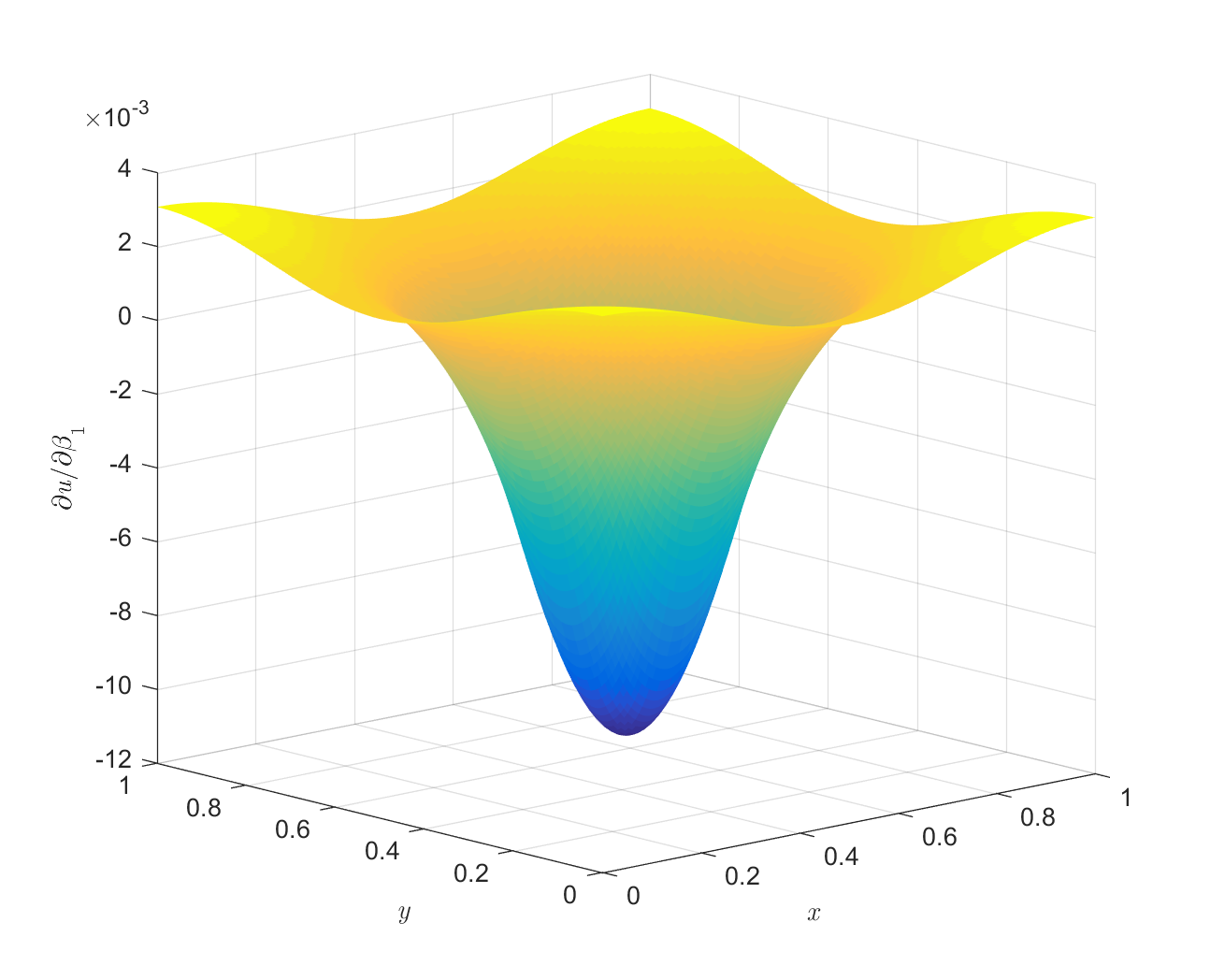}
  \caption{Finite element approximation of $\partial_{\beta_1}u$ under a refined mesh $h=1/128$ in the membrane problem}\label{derivdisp1}
\end{figure}

Under a refined mesh $h=1/128$, one has the solution of $\partial_{\beta_1}u$ as shown in Figure \ref{derivdisp1}, and the reference value of $J_1$ is found to be $-8.762\times 10^{-3}$.

\begin{figure}[hbt]
  \begin{minipage}[t]{0.48\linewidth}
  \centering
   \includegraphics[width=\textwidth]{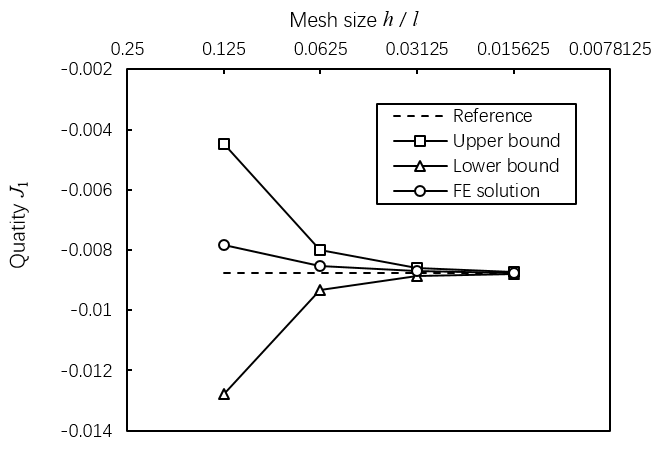}
  \caption{FE solutions, upper and lower bounds of $J_1$ ($\xi=1.0$) for the membrane}\label{memqoi1}
  \end{minipage}
  \hfill
  \begin{minipage}[t]{0.48\linewidth}
  \centering
   \includegraphics[width=\textwidth]{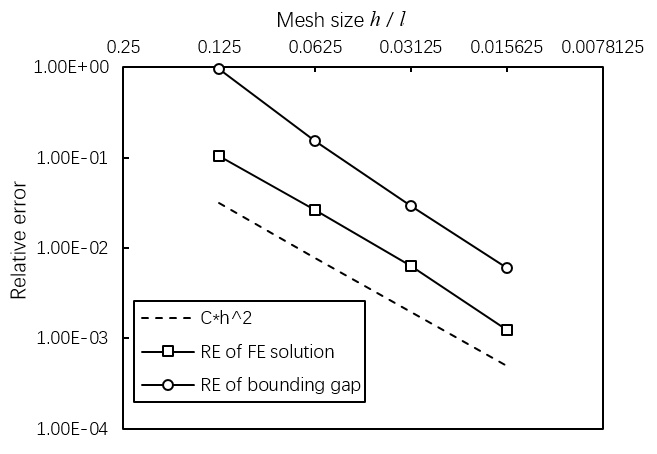}
  \caption{Relative errors of FE solutions and bounding gaps of $J_1$ ($\xi=1.0$) for the membrane}\label{memRE1}
  \end{minipage}
\end{figure}

In this case, the range of $\xi$ is $0<\xi<(2/\bar{\beta}_1)=2.0$. Taking $\xi=1.0$, it is seen in Figure \ref{memqoi1} that strict bounding property of $J_1$ is achieved under meshes with various sizes. The relative errors of FE solutions and bounding gaps of the quantity are shown in Figure \ref{memRE1}, in which both the finite element solution and the bounding gap have the same convergence rate of $h^2$, implying the super-convergent asymptotic property of the proposed error bounds.

\textbf{Case 2:} $J_2=\partial_{\beta_2}\left[\left(\int_{\Omega_{QoI}}u\right)/\left|\Omega_{QoI}\right|\right]=64\int_{\Omega_{QoI}}\partial_{\beta_2}u$

\begin{figure}[htb]
  \centering
  \includegraphics[scale=0.3]{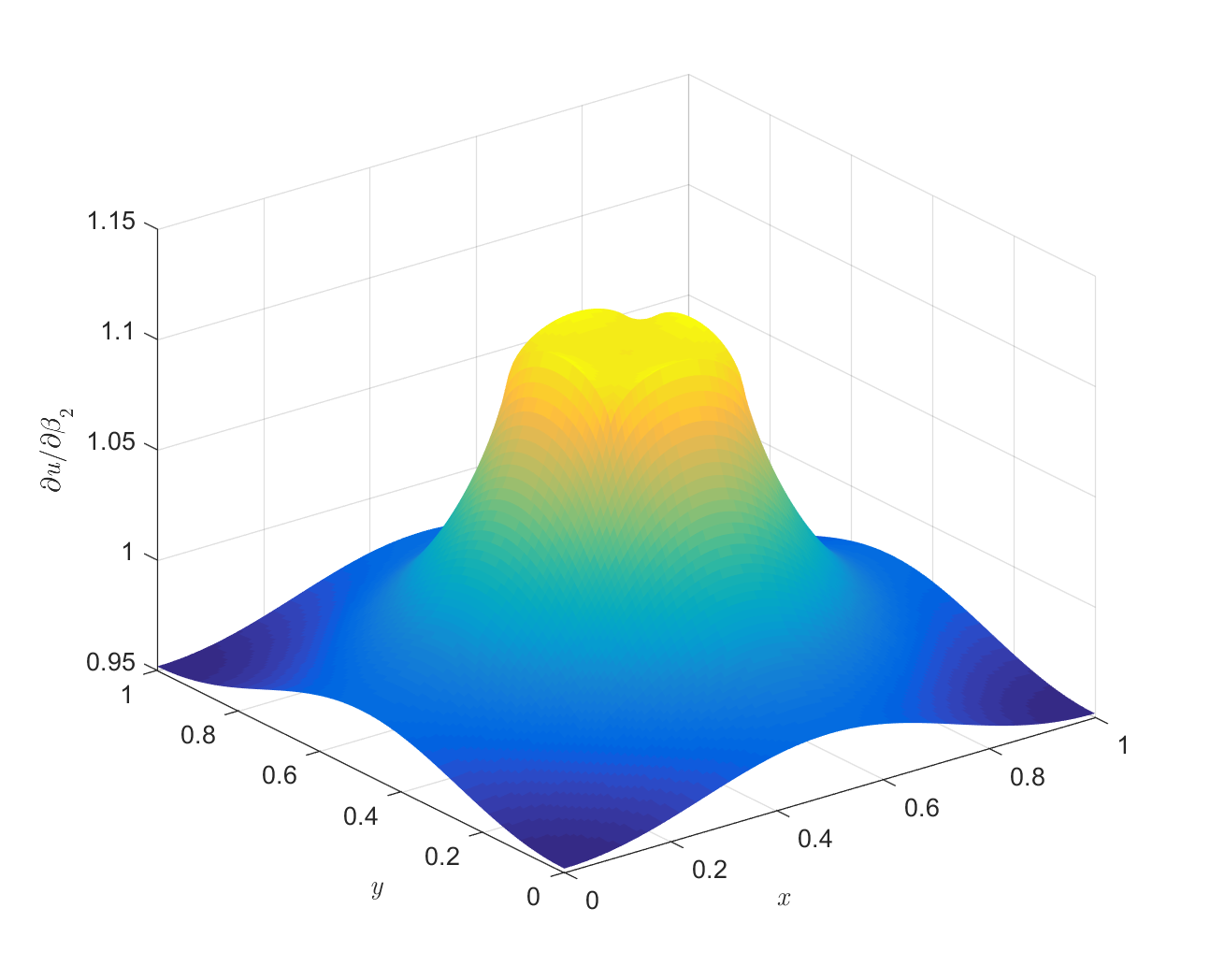}
  \caption{Finite element approximation of $\partial_{\beta_2}u$ under a refined mesh $h=1/128$ in the membrane problem}\label{derivdisp2}
\end{figure}

Under a refined mesh $h=1/128$, one has the solution of $\partial_{\beta_2}u$ as shown in Figure \ref{derivdisp2}, and the reference value of $J_2$ is found to be $1.1060$.

\begin{figure}[hbt]
  \begin{minipage}[t]{0.48\linewidth}
  \centering
   \includegraphics[width=\textwidth]{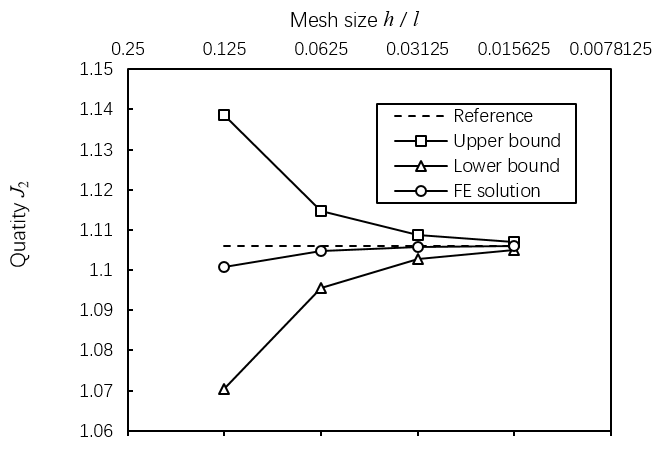}
  \caption{FE solutions, upper and lower bounds of $J_2$ for the membrane}\label{memqoi2}
  \end{minipage}
  \hfill
  \begin{minipage}[t]{0.48\linewidth}
  \centering
   \includegraphics[width=\textwidth]{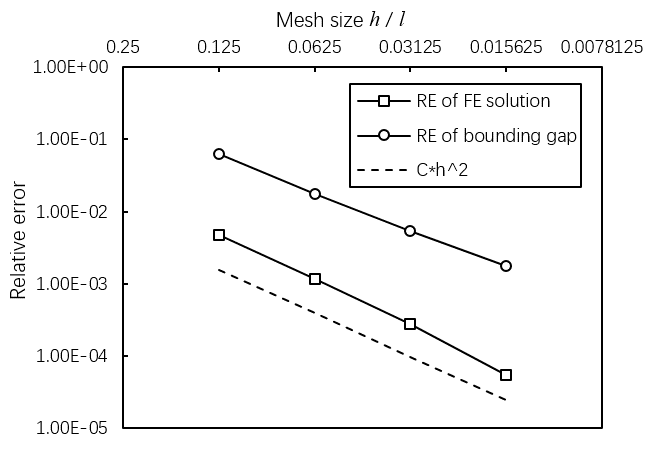}
  \caption{Relative errors of FE solutions and bounding gaps of $J_2$ for the membrane}\label{memRE2}
  \end{minipage}
\end{figure}

Since the load is exclusively parameterized in this case, $\xi$ is irrelevant. The derivative term $\langle \partial_{\beta_2}f,\cdot\rangle$ can be explicitly written as
\begin{equation*}
\langle \partial_{\beta_2}f,v\rangle=\int_{\partial \Omega_{QoI}}1\cdot v\,,\quad v\in \mathcal{V}\,.
\end{equation*}
Numerical results for quantity $J_2$, including the FE solutions, upper and lower bounds under different mesh densities are illustrated in Figure \ref{memqoi2}, and the strict bounding property of the proposed estimation technique for sensitivity derivative is displayed again. Figure \ref{memRE2} shows the relative errors of FE solutions and bounding gaps of this quantity.

\section{Conclusions}

In the sense of dual variational principles, the CRE-based goal-oriented error estimation has been extended to the cases with non-symmetric bilinear forms, and applied to static response sensitivity analysis of linear structural problems. Strict upper and lower bounds of quantities with respect to the sensitivity derivative fields, including the sensitivity derivatives of some structural response quantities, have been acquired by the proposed technique. The present goal-oriented error estimation is employed in sensitivity analysis of a Bernoulli-Euler beam problem and a membrane problem, and numerical results have validated the strict bounding property and the same convergence rate of the proposed bounds as the quantity error.

\section*{Acknowledgment}

This work was supported by the National Natural Science Foundation of China under Grant No. 51378294.

\end{document}